\documentstyle[11pt]{article}

\title{\Large{ \bf A categorical construction of  2-dimensional extended
           Topological Quantum Field Theory}}

\author{Vishvajit V S Gautam\thanks{E-mail : vishvajit@imsc.res.in; gautamvvs@yahoo.com }\\
\small \sl  The Institute of Mathematical Sciences\\
\small \sl  CIT Campus Taramani, Chennai - 600113 INDIA}
\date{ }
\begin{document}
\newtheorem{defi}{\bf Definition}[section] 
\newtheorem{th}[defi]{\bf Theorem}
\newtheorem{pro}[defi]{\bf Proposition}
\newtheorem{lemma}[defi]{\bf Lemma}
\newtheorem{coro}[defi]{\bf Corollary}
\newtheorem{ram}[defi]{\bf Remark}
\newtheorem{pn}[subsubsection]{\bf Proposition}
\newtheorem{zm}[subsubsection]{\bf Theorem}
\newtheorem{cy}[subsubsection]{\bf Corollary}
\newtheorem{df}[subsubsection]{\bf Definition}
\newtheorem{rk}[subsubsection]{\bf Remark}
\maketitle

\begin{abstract}

In this  paper we propose a naive construction of 2-dimensional
extended topological quantum field theories (TQFTs), which can be further
generalized to the higher-dimension extended TQFTs.\\



\end{abstract}

\section{Introduction}

The notion of TQFT is related with the study of path integral for
Lagrangian with topological invariance. In some sense a 
TQFT gives a topological invariant to $n$-dimensional
manifolds. For example 3-dimensional TQFT assigns a 3-manifold $M$ to
a numerical invariant ${\tau}(M)$ (which may be a complex number)
 such that if ${\tau}(M) \neq
{\tau}{(M')}$ for 3-manifolds $M$ and $M'$ then $M$ and $M'$ are not diffeomorphic. 
A 3-dimensional TQFT also gives numerical invariants of knots, links
and ribbon graphs. This theory was introduced by Witten in 1988 to
describe a class of quantum field theory whose action is
diffeomorphism invariant such as $\tau (M)$  which not
only assigns a numerical value to manifold but also preserves the
embedded ribbon graphs or any other structure defined over the manifold.  Atiyah in 1988
formulated an axiomatic setup for TQFTs. Independently and at
about the same time G. Segal formulates a mathematical definition of
{\em conformal field theories} (CFTs), which is very similarly
based on categories and functors.

Categories play a centre role in mathematical formulation of TQFT. An $n$-dimensional TQFT is defined as a {\em
  monoidal functor} from the category of oriented {\em $n$-cobordisms}
  with disjoint union as tensor product to the category {\bf Vect} of
  finite dimensional vector spaces with usual tensor product of vector spaces.
Any modification in cobordism category may leads to a modification in 
TQFT. This modification can be think as a {\em extended} version of TQFT.
For example in Chern-Simons Witten TQFT cobordisms are supplied with
  some extra structures [4].

Similar to the case of TQFTs  there are many ways to define
{\em extended TQFTs}. Kerler and Lyubashenko in [8] introduce a
notion of {\em extended TQFTs} which involves higher category theory, namely double
categories and double functors. Involvement of higher-dimensional
algebra indicate that extended TQFTs will preserve extra
information. Due to double functority, their extended TQFTs
contains both Atiyah notion of a TQFT in dimension three and Segal's notion of
CFT as special cases, though they appear on different
categorical levels.

The role of higher-dimensional algebra is clear from the various
constructions of extended TQFTs. Baez and Dolan in
[2] outline a program in
which $n$-dimensional TQFTs are described as {\em $n$-category}
representation. They described an $n$-dimensional extended TQFT as
a weak $n$-functor from the free stable weak $n$-category with duals
of one objects to {\bf nHilb} the category of $n$-Hilbert spaces,
which preserve  all levels of duality.

This paper gives a naive categorical construction of a 2-dimensional
extended TQFT, which is different from the known constructions of
$n$-dimensional extended TQFT for $n = 2$. We use the notion of
internal categories {\em 2-Vector spaces} which is in agreement with the
2-Vector spaces defined  by Baez and Crans in [1].

\section{2- Categories and  semistrict monoidal 2-categories}

2-categories  are the first prototype of higher-dimensional algebra. Importance of symmetric monoidal 2-categories,
braided monoidal 2-categories is evident from the recent developments
in higher-dimensional algebra, see  Baez and Neuchl [3], Day and
Street [5], and Kapranov and Voevodsky [7], and reference therein. Semistrict
monoidal 2-categories are the base categories for our description
of an extended TQFTs. Instead of using the weak version of
monoidal 2-categories we will use the semistrict version of monoidal
2-categories, because they are now better understood, and a
coherence theorem (cf. [6]) for weak categories says they are equivalent to
semistrict ones.

\subsection{2-Category}

A 2-category $\cal C$ consists of following data:\\

$\bullet$ a class ${\cal C}_{0}$ of objects $A$, $B$, $C$, ... which are
  called {\em $0$-cells};

The objects or 0-cells and arrows or 1-cells form a category, called
the underlying category of $\cal C$ which we also
denote by $\cal C$, with identities ${1_{A}} : A \longrightarrow A$.
\\

$\bullet$ for each pair $A$,$B$ in ${\cal C}_{0}$, a small category ${\cal
  C}_{{1}_{(A,B)}}$ whose objects or {\em $1$-cells} are morphisms $f: A
  \longrightarrow B$ etc., and arrows or {\em $2$-cells} are morphisms of morphisms from
  $A$ to $B$, which we denote by $\alpha$, $\beta$, $\gamma$...;\\
  A $2$-cell pictured as\\
\\

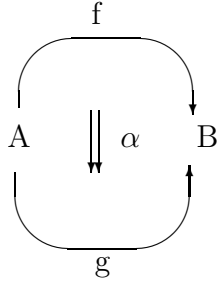
\setlength{\unitlength}{2000sp}%
\begingroup\makeatletter\ifx\SetFigFont\undefined
\def\x#1#2#3#4#5#6#7\relax{\def\x{#1#2#3#4#5#6}}%
\expandafter\x\fmtname xxxxxx\relax \def\y{splain}%
\ifx\x\y   
\gdef\SetFigFont#1#2#3{%
  \ifnum #1<17\tiny\else \ifnum #1<20\small\else
  \ifnum #1<24\normalsize\else \ifnum #1<29\large\else
  \ifnum #1<34\Large\else \ifnum #1<41\LARGE\else
     \huge\fi\fi\fi\fi\fi\fi
  \csname #3\endcsname}%
\else
\gdef\SetFigFont#1#2#3{\begingroup
  \count@#1\relax \ifnum 25<\count@\count@25\fi
  \def\x{\endgroup\@setsize\SetFigFont{#2pt}}%
  \expandafter\x
    \csname \romannumeral\the\count@ pt\expandafter\endcsname
    \csname @\romannumeral\the\count@ pt\endcsname
  \csname #3\endcsname}%
\fi
\fi\endgroup
\begin{picture}(2340,3235)(1351,-2666)
\thinlines
\put(3646,-691){\vector( 0,-1){0}}
\put(2558,-691){\oval(2176,1896)[tr]}
\put(2558,-601){\oval(2144,1716)[tl]}
\put(2517,-1411){\oval(2152,1890)[bl]}
\put(2517,-1321){\oval(2168,2070)[br]}
\put(3601,-1321){\vector( 0, 1){0}}
\put(2386,-646){\vector( 0,-1){765}}
\put(2476,-646){\vector( 0,-1){765}}
\put(1351,-1096){\makebox(0,0)[lb]{\smash{\SetFigFont{12}{14.4}{rm}A }}}
\put(3691,-1096){\makebox(0,0)[lb]{\smash{\SetFigFont{12}{14.4}{rm}B}}}
\put(2386,434){\makebox(0,0)[lb]{\smash{\SetFigFont{12}{14.4}{rm}f}}}
\put(2431,-2626){\makebox(0,0)[lb]{\smash{\SetFigFont{12}{14.4}{rm}g}}}
\put(2746,-1096){\makebox(0,0)[lb]{\smash{\SetFigFont{12}{14.4}{rm}$\alpha$}}}
\end{picture}

For any pair 2-cells $\alpha$ , $\beta $ in ${\cal C}_{{1}_{(A,B)}}$ the 2-cells
composition under which ${\cal C}_{{1}_{(A,B)}}$ form a category  is
called {\em vertical composition}.\\
A 2-cells vertical composition as displayed in:\\

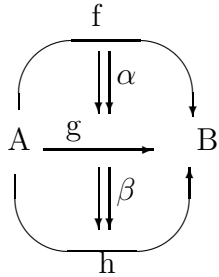
\setlength{\unitlength}{2000sp}%
\begingroup\makeatletter\ifx\SetFigFont\undefined
\def\x#1#2#3#4#5#6#7\relax{\def\x{#1#2#3#4#5#6}}%
\expandafter\x\fmtname xxxxxx\relax \def\y{splain}%
\ifx\x\y   
\gdef\SetFigFont#1#2#3{%
  \ifnum #1<17\tiny\else \ifnum #1<20\small\else
  \ifnum #1<24\normalsize\else \ifnum #1<29\large\else
  \ifnum #1<34\Large\else \ifnum #1<41\LARGE\else
     \huge\fi\fi\fi\fi\fi\fi
  \csname #3\endcsname}%
\else
\gdef\SetFigFont#1#2#3{\begingroup
  \count@#1\relax \ifnum 25<\count@\count@25\fi
  \def\x{\endgroup\@setsize\SetFigFont{#2pt}}%
  \expandafter\x
    \csname \romannumeral\the\count@ pt\expandafter\endcsname
    \csname @\romannumeral\the\count@ pt\endcsname
  \csname #3\endcsname}%
\fi
\fi\endgroup
\begin{picture}(2340,3195)(1351,-2626)
\thinlines
\put(3646,-691){\vector( 0,-1){0}}
\put(2558,-691){\oval(2176,1896)[tr]}
\put(2558,-601){\oval(2144,1716)[tl]}
\put(2517,-1411){\oval(2152,1890)[bl]}
\put(2517,-1321){\oval(2168,2070)[br]}
\put(3601,-1321){\vector( 0, 1){0}}
\put(2611,119){\vector( 0,-1){765}}
\put(2476,119){\vector( 0,-1){765}}
\put(2611,-1321){\vector( 0,-1){765}}
\put(2476,-1321){\vector( 0,-1){765}}
\put(1801,-1096){\vector( 1, 0){1350}}
\put(1351,-1096){\makebox(0,0)[lb]{\smash{\SetFigFont{12}{14.4}{rm}A }}}
\put(3691,-1096){\makebox(0,0)[lb]{\smash{\SetFigFont{12}{14.4}{rm}B}}}
\put(2386,434){\makebox(0,0)[lb]{\smash{\SetFigFont{12}{14.4}{rm}f}}}
\put(2701,-286){\makebox(0,0)[lb]{\smash{\SetFigFont{12}{14.4}{rm}$\alpha$}}}
\put(2701,-1726){\makebox(0,0)[lb]{\smash{\SetFigFont{12}{14.4}{rm}$\beta$}}}
\put(2476,-2626){\makebox(0,0)[lb]{\smash{\SetFigFont{12}{14.4}{rm}h}}}
\put(2071,-916){\makebox(0,0)[lb]{\smash{\SetFigFont{12}{14.4}{rm}g}}}
\end{picture}

The vertical composite $f \Rightarrow h$ is denoted by $\beta \circ
\alpha$. Identities of  ${\cal C}_{{1}_{(A,B)}}$ are denoted by\\

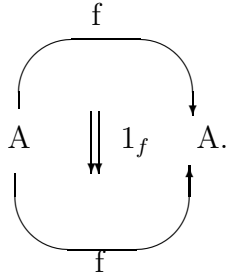
\setlength{\unitlength}{2000sp}%
\begingroup\makeatletter\ifx\SetFigFont\undefined
\def\x#1#2#3#4#5#6#7\relax{\def\x{#1#2#3#4#5#6}}%
\expandafter\x\fmtname xxxxxx\relax \def\y{splain}%
\ifx\x\y   
\gdef\SetFigFont#1#2#3{%
  \ifnum #1<17\tiny\else \ifnum #1<20\small\else
  \ifnum #1<24\normalsize\else \ifnum #1<29\large\else
  \ifnum #1<34\Large\else \ifnum #1<41\LARGE\else
     \huge\fi\fi\fi\fi\fi\fi
  \csname #3\endcsname}%
\else
\gdef\SetFigFont#1#2#3{\begingroup
  \count@#1\relax \ifnum 25<\count@\count@25\fi
  \def\x{\endgroup\@setsize\SetFigFont{#2pt}}%
  \expandafter\x
    \csname \romannumeral\the\count@ pt\expandafter\endcsname
    \csname @\romannumeral\the\count@ pt\endcsname
  \csname #3\endcsname}%
\fi
\fi\endgroup
\begin{picture}(2340,3235)(1351,-2666)
\thinlines
\put(3646,-691){\vector( 0,-1){0}}
\put(2558,-691){\oval(2176,1896)[tr]}
\put(2558,-601){\oval(2144,1716)[tl]}
\put(2517,-1411){\oval(2152,1890)[bl]}
\put(2517,-1321){\oval(2168,2070)[br]}
\put(3601,-1321){\vector( 0, 1){0}}
\put(2386,-646){\vector( 0,-1){765}}
\put(2476,-646){\vector( 0,-1){765}}
\put(1351,-1096){\makebox(0,0)[lb]{\smash{\SetFigFont{12}{14.4}{rm}A }}}
\put(3691,-1096){\makebox(0,0)[lb]{\smash{\SetFigFont{12}{14.4}{rm}A.}}}
\put(2386,434){\makebox(0,0)[lb]{\smash{\SetFigFont{12}{14.4}{rm}f}}}
\put(2431,-2626){\makebox(0,0)[lb]{\smash{\SetFigFont{12}{14.4}{rm}f}}}
\put(2746,-1096){\makebox(0,0)[lb]{\smash{\SetFigFont{12}{14.4}{rm}$1_{f}$}}}
\end{picture}

Between any pair of  0-cells $A$ and $B$ there are  2-cells $\alpha$
etc. We can compose them  under  another 2-cells operation known as
{\em horizontal composition}. In this case we have
\\

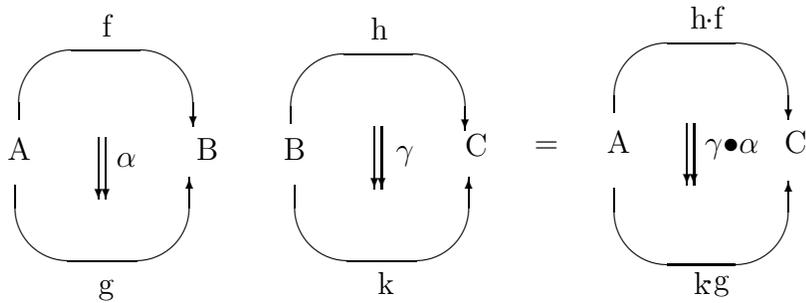
\setlength{\unitlength}{2000sp}%
\begingroup\makeatletter\ifx\SetFigFont\undefined
\def\x#1#2#3#4#5#6#7\relax{\def\x{#1#2#3#4#5#6}}%
\expandafter\x\fmtname xxxxxx\relax \def\y{splain}%
\ifx\x\y   
\gdef\SetFigFont#1#2#3{%
  \ifnum #1<17\tiny\else \ifnum #1<20\small\else
  \ifnum #1<24\normalsize\else \ifnum #1<29\large\else
  \ifnum #1<34\Large\else \ifnum #1<41\LARGE\else
     \huge\fi\fi\fi\fi\fi\fi
  \csname #3\endcsname}%
\else
\gdef\SetFigFont#1#2#3{\begingroup
  \count@#1\relax \ifnum 25<\count@\count@25\fi
  \def\x{\endgroup\@setsize\SetFigFont{#2pt}}%
  \expandafter\x
    \csname \romannumeral\the\count@ pt\expandafter\endcsname
    \csname @\romannumeral\the\count@ pt\endcsname
  \csname #3\endcsname}%
\fi
\fi\endgroup
\begin{picture}(9697,3460)(1351,-2801)
\thinlines
\put(3646,-691){\vector( 0,-1){0}}
\put(2558,-691){\oval(2176,1896)[tr]}
\put(2558,-601){\oval(2144,1716)[tl]}
\put(2517,-1411){\oval(2152,1890)[bl]}
\put(2517,-1321){\oval(2168,2070)[br]}
\put(3601,-1321){\vector( 0, 1){0}}
\put(7021,-736){\vector( 0,-1){0}}
\put(5933,-736){\oval(2176,1896)[tr]}
\put(5933,-646){\oval(2144,1716)[tl]}
\put(5982,-1411){\oval(2152,1890)[bl]}
\put(5982,-1321){\oval(2168,2070)[br]}
\put(7066,-1321){\vector( 0, 1){0}}
\put(11026,-601){\vector( 0,-1){0}}
\put(9938,-601){\oval(2176,1896)[tr]}
\put(9938,-511){\oval(2144,1716)[tl]}
\put(9942,-1456){\oval(2152,1890)[bl]}
\put(9942,-1366){\oval(2168,2070)[br]}
\put(11026,-1366){\vector( 0, 1){0}}
\put(5896,-691){\vector( 0,-1){765}}
\put(5986,-691){\vector( 0,-1){765}}
\put(2476,-826){\vector( 0,-1){765}}
\put(2566,-826){\vector( 0,-1){765}}
\put(9856,-646){\vector( 0,-1){765}}
\put(9766,-646){\vector( 0,-1){765}}
\put(1351,-1096){\makebox(0,0)[lb]{\smash{\SetFigFont{12}{14.4}{rm}A }}}
\put(3691,-1096){\makebox(0,0)[lb]{\smash{\SetFigFont{12}{14.4}{rm}B}}}
\put(4771,-1096){\makebox(0,0)[lb]{\smash{\SetFigFont{12}{14.4}{rm}B}}}
\put(7021,-1051){\makebox(0,0)[lb]{\smash{\SetFigFont{12}{14.4}{rm}C}}}
\put(6166,-1096){\makebox(0,0)[lb]{\smash{\SetFigFont{12}{14.4}{rm}$\gamma $}}}
\put(2701,-1186){\makebox(0,0)[lb]{\smash{\SetFigFont{12}{14.4}{rm}$\alpha$}}}
\put(2521,434){\makebox(0,0)[lb]{\smash{\SetFigFont{12}{14.4}{rm}f}}}
\put(5851,389){\makebox(0,0)[lb]{\smash{\SetFigFont{12}{14.4}{rm}h}}}
\put(5941,-2761){\makebox(0,0)[lb]{\smash{\SetFigFont{12}{14.4}{rm}k}}}
\put(2476,-2761){\makebox(0,0)[lb]{\smash{\SetFigFont{12}{14.4}{rm}g}}}
\put(7876,-1006){\makebox(0,0)[lb]{\smash{\SetFigFont{12}{14.4}{rm}$=$ }}}
\put(8776,-1006){\makebox(0,0)[lb]{\smash{\SetFigFont{12}{14.4}{rm}A }}}
\put(9991,-1006){\makebox(0,0)[lb]{\smash{\SetFigFont{12}{14.4}{rm}$\gamma$}}}
\put(10216,-1006){\makebox(0,0)[lb]{\smash{\SetFigFont{12}{14.4}{rm}$\bullet
      \alpha$}}}
\put(9811,524){\makebox(0,0)[lb]{\smash{\SetFigFont{12}{14.4}{rm}h}}}
\put(9991,524){\makebox(0,0)[lb]{\smash{\SetFigFont{12}{14.4}{rm}$\cdot$f}}}
\put(9991,-2716){\makebox(0,0)[lb]{\smash{\SetFigFont{12}{14.4}{rm}$\cdot$g}}}
\put(9856,-2761){\makebox(0,0)[lb]{\smash{\SetFigFont{12}{14.4}{rm}k}}}
\put(10981,-1006){\makebox(0,0)[lb]{\smash{\SetFigFont{12}{14.4}{rm}C}}}
\end{picture}

Under this composite $\gamma \bullet \alpha : h \cdot f \Rightarrow k
\cdot g $ law the 2-cells form a category, with identities\\

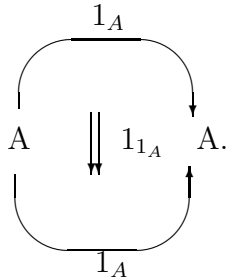
\setlength{\unitlength}{2000sp}%
\begingroup\makeatletter\ifx\SetFigFont\undefined
\def\x#1#2#3#4#5#6#7\relax{\def\x{#1#2#3#4#5#6}}%
\expandafter\x\fmtname xxxxxx\relax \def\y{splain}%
\ifx\x\y   
\gdef\SetFigFont#1#2#3{%
  \ifnum #1<17\tiny\else \ifnum #1<20\small\else
  \ifnum #1<24\normalsize\else \ifnum #1<29\large\else
  \ifnum #1<34\Large\else \ifnum #1<41\LARGE\else
     \huge\fi\fi\fi\fi\fi\fi
  \csname #3\endcsname}%
\else
\gdef\SetFigFont#1#2#3{\begingroup
  \count@#1\relax \ifnum 25<\count@\count@25\fi
  \def\x{\endgroup\@setsize\SetFigFont{#2pt}}%
  \expandafter\x
    \csname \romannumeral\the\count@ pt\expandafter\endcsname
    \csname @\romannumeral\the\count@ pt\endcsname
  \csname #3\endcsname}%
\fi
\fi\endgroup
\begin{picture}(2340,3235)(1351,-2666)
\thinlines
\put(3646,-691){\vector( 0,-1){0}}
\put(2558,-691){\oval(2176,1896)[tr]}
\put(2558,-601){\oval(2144,1716)[tl]}
\put(2517,-1411){\oval(2152,1890)[bl]}
\put(2517,-1321){\oval(2168,2070)[br]}
\put(3601,-1321){\vector( 0, 1){0}}
\put(2386,-646){\vector( 0,-1){765}}
\put(2476,-646){\vector( 0,-1){765}}
\put(1351,-1096){\makebox(0,0)[lb]{\smash{\SetFigFont{12}{14.4}{rm}A }}}
\put(3691,-1096){\makebox(0,0)[lb]{\smash{\SetFigFont{12}{14.4}{rm}A.}}}
\put(2386,434){\makebox(0,0)[lb]{\smash{\SetFigFont{12}{14.4}{rm}$1_{A}$}}}
\put(2431,-2626){\makebox(0,0)[lb]{\smash{\SetFigFont{12}{14.4}{rm}$1_{A}$}}}
\put(2746,-1096){\makebox(0,0)[lb]{\smash{\SetFigFont{12}{14.4}{rm}$1_{1_{A}}$}}}
\end{picture}

Finally the horizontal and the vertical compositions are related with
following conditions: \\
\\

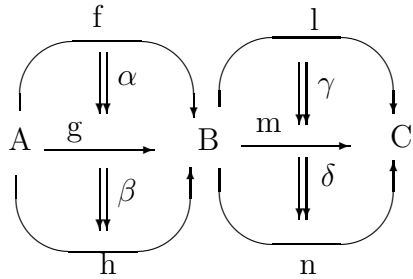
\setlength{\unitlength}{2000sp}%
\begingroup\makeatletter\ifx\SetFigFont\undefined
\def\x#1#2#3#4#5#6#7\relax{\def\x{#1#2#3#4#5#6}}%
\expandafter\x\fmtname xxxxxx\relax \def\y{splain}%
\ifx\x\y   
\gdef\SetFigFont#1#2#3{%
  \ifnum #1<17\tiny\else \ifnum #1<20\small\else
  \ifnum #1<24\normalsize\else \ifnum #1<29\large\else
  \ifnum #1<34\Large\else \ifnum #1<41\LARGE\else
     \huge\fi\fi\fi\fi\fi\fi
  \csname #3\endcsname}%
\else
\gdef\SetFigFont#1#2#3{\begingroup
  \count@#1\relax \ifnum 25<\count@\count@25\fi
  \def\x{\endgroup\@setsize\SetFigFont{#2pt}}%
  \expandafter\x
    \csname \romannumeral\the\count@ pt\expandafter\endcsname
    \csname @\romannumeral\the\count@ pt\endcsname
  \csname #3\endcsname}%
\fi
\fi\endgroup
\begin{picture}(4792,3195)(1351,-2626)
\thinlines
\put(3646,-691){\vector( 0,-1){0}}
\put(2558,-691){\oval(2176,1896)[tr]}
\put(2558,-601){\oval(2144,1716)[tl]}
\put(2517,-1411){\oval(2152,1890)[bl]}
\put(2517,-1321){\oval(2168,2070)[br]}
\put(3601,-1321){\vector( 0, 1){0}}
\put(6121,-646){\vector( 0,-1){0}}
\put(5033,-646){\oval(2176,1896)[tr]}
\put(5033,-556){\oval(2144,1716)[tl]}
\put(5037,-1321){\oval(2152,1890)[bl]}
\put(5037,-1231){\oval(2168,2070)[br]}
\put(6121,-1231){\vector( 0, 1){0}}
\put(2476,119){\vector( 0,-1){765}}
\put(2476,-1321){\vector( 0,-1){765}}
\put(1801,-1096){\vector( 1, 0){1350}}
\put(4231,-1051){\vector( 1, 0){1350}}
\put(4951,-16){\vector( 0,-1){765}}
\put(4951,-1186){\vector( 0,-1){765}}
\put(5041,-1186){\vector( 0,-1){765}}
\put(5041,-16){\vector( 0,-1){765}}
\put(2566,119){\vector( 0,-1){765}}
\put(2566,-1321){\vector( 0,-1){765}}
\put(1351,-1096){\makebox(0,0)[lb]{\smash{\SetFigFont{12}{14.4}{rm}A }}}
\put(3691,-1096){\makebox(0,0)[lb]{\smash{\SetFigFont{12}{14.4}{rm}B}}}
\put(2386,434){\makebox(0,0)[lb]{\smash{\SetFigFont{12}{14.4}{rm}f}}}
\put(2701,-286){\makebox(0,0)[lb]{\smash{\SetFigFont{12}{14.4}{rm}$\alpha$}}}
\put(2701,-1726){\makebox(0,0)[lb]{\smash{\SetFigFont{12}{14.4}{rm}$\beta$}}}
\put(2476,-2626){\makebox(0,0)[lb]{\smash{\SetFigFont{12}{14.4}{rm}h}}}
\put(2071,-916){\makebox(0,0)[lb]{\smash{\SetFigFont{12}{14.4}{rm}g}}}
\put(6076,-1051){\makebox(0,0)[lb]{\smash{\SetFigFont{12}{14.4}{rm}C}}}
\put(4951,-2626){\makebox(0,0)[lb]{\smash{\SetFigFont{12}{14.4}{rm}n}}}
\put(5086,389){\makebox(0,0)[lb]{\smash{\SetFigFont{12}{14.4}{rm}l}}}
\put(5176,-376){\makebox(0,0)[lb]{\smash{\SetFigFont{12}{14.4}{rm}$\gamma$}}}
\put(5221,-1546){\makebox(0,0)[lb]{\smash{\SetFigFont{12}{14.4}{rm}$\delta$}}}
\put(4411,-871){\makebox(0,0)[lb]{\smash{\SetFigFont{12}{14.4}{rm}m}}}
\end{picture}
\\
\[ (\delta \bullet \beta ) \circ (\gamma \bullet \alpha) = (\delta
\circ \gamma ) \bullet (\beta \circ \alpha) \]

In the situation\\

\setlength{\unitlength}{2000sp}%
\begingroup\makeatletter\ifx\SetFigFont\undefined
\def\x#1#2#3#4#5#6#7\relax{\def\x{#1#2#3#4#5#6}}%
\expandafter\x\fmtname xxxxxx\relax \def\y{splain}%
\ifx\x\y   
\gdef\SetFigFont#1#2#3{%
  \ifnum #1<17\tiny\else \ifnum #1<20\small\else
  \ifnum #1<24\normalsize\else \ifnum #1<29\large\else
  \ifnum #1<34\Large\else \ifnum #1<41\LARGE\else
     \huge\fi\fi\fi\fi\fi\fi
  \csname #3\endcsname}%
\else
\gdef\SetFigFont#1#2#3{\begingroup
  \count@#1\relax \ifnum 25<\count@\count@25\fi
  \def\x{\endgroup\@setsize\SetFigFont{#2pt}}%
  \expandafter\x
    \csname \romannumeral\the\count@ pt\expandafter\endcsname
    \csname @\romannumeral\the\count@ pt\endcsname
  \csname #3\endcsname}%
\fi
\fi\endgroup
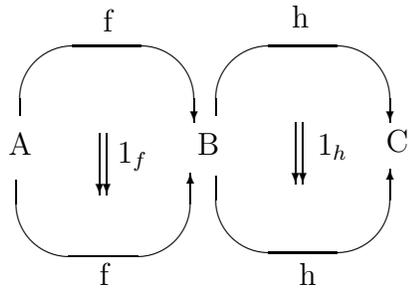
\begin{picture}(4747,3330)(1351,-2716)
\thinlines
\put(3646,-691){\vector( 0,-1){0}}
\put(2558,-691){\oval(2176,1896)[tr]}
\put(2558,-601){\oval(2144,1716)[tl]}
\put(2517,-1411){\oval(2152,1890)[bl]}
\put(2517,-1321){\oval(2168,2070)[br]}
\put(3601,-1321){\vector( 0, 1){0}}
\put(4992,-1366){\oval(2152,1890)[bl]}
\put(4992,-1276){\oval(2168,2070)[br]}
\put(6076,-1276){\vector( 0, 1){0}}
\put(6076,-691){\vector( 0,-1){0}}
\put(4988,-691){\oval(2176,1896)[tr]}
\put(4988,-601){\oval(2144,1716)[tl]}
\put(2476,-826){\vector( 0,-1){765}}
\put(2566,-826){\vector( 0,-1){765}}
\put(4996,-691){\vector( 0,-1){765}}
\put(4906,-691){\vector( 0,-1){765}}
\put(1351,-1096){\makebox(0,0)[lb]{\smash{\SetFigFont{12}{14.4}{rm}A }}}
\put(3691,-1096){\makebox(0,0)[lb]{\smash{\SetFigFont{12}{14.4}{rm}B}}}
\put(2701,-1186){\makebox(0,0)[lb]{\smash{\SetFigFont{12}{14.4}{rm}$1_{f}$}}}
\put(2521,434){\makebox(0,0)[lb]{\smash{\SetFigFont{12}{14.4}{rm}f}}}
\put(4861,479){\makebox(0,0)[lb]{\smash{\SetFigFont{12}{14.4}{rm}h}}}
\put(5176,-1096){\makebox(0,0)[lb]{\smash{\SetFigFont{12}{14.4}{rm}$1_{h}$}}}
\put(6031,-1051){\makebox(0,0)[lb]{\smash{\SetFigFont{12}{14.4}{rm}C}}}
\put(2476,-2716){\makebox(0,0)[lb]{\smash{\SetFigFont{12}{14.4}{rm}f}}}
\put(4951,-2716){\makebox(0,0)[lb]{\smash{\SetFigFont{12}{14.4}{rm}h}}}
\end{picture}
\\
 we require that horizontal composite of two vertical identities is
 itself a vertical identity i.e. $1_{h} \bullet 1_{f} =  1_{h \cdot f}$.

This structure also provides a horizontal composite of a 2-cell with
1-cell\\

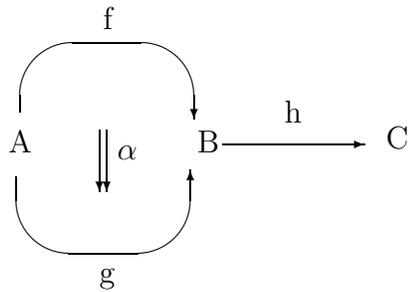
\setlength{\unitlength}{2000sp}%
\begingroup\makeatletter\ifx\SetFigFont\undefined
\def\x#1#2#3#4#5#6#7\relax{\def\x{#1#2#3#4#5#6}}%
\expandafter\x\fmtname xxxxxx\relax \def\y{splain}%
\ifx\x\y   
\gdef\SetFigFont#1#2#3{%
  \ifnum #1<17\tiny\else \ifnum #1<20\small\else
  \ifnum #1<24\normalsize\else \ifnum #1<29\large\else
  \ifnum #1<34\Large\else \ifnum #1<41\LARGE\else
     \huge\fi\fi\fi\fi\fi\fi
  \csname #3\endcsname}%
\else
\gdef\SetFigFont#1#2#3{\begingroup
  \count@#1\relax \ifnum 25<\count@\count@25\fi
  \def\x{\endgroup\@setsize\SetFigFont{#2pt}}%
  \expandafter\x
    \csname \romannumeral\the\count@ pt\expandafter\endcsname
    \csname @\romannumeral\the\count@ pt\endcsname
  \csname #3\endcsname}%
\fi
\fi\endgroup
\begin{picture}(4680,3325)(1351,-2756)
\thinlines
\put(3646,-691){\vector( 0,-1){0}}
\put(2558,-691){\oval(2176,1896)[tr]}
\put(2558,-601){\oval(2144,1716)[tl]}
\put(2517,-1411){\oval(2152,1890)[bl]}
\put(2517,-1321){\oval(2168,2070)[br]}
\put(3601,-1321){\vector( 0, 1){0}}
\put(2476,-826){\vector( 0,-1){765}}
\put(2566,-826){\vector( 0,-1){765}}
\put(4006,-1006){\vector( 1, 0){1755}}
\put(1351,-1096){\makebox(0,0)[lb]{\smash{\SetFigFont{12}{14.4}{rm}A }}}
\put(3691,-1096){\makebox(0,0)[lb]{\smash{\SetFigFont{12}{14.4}{rm}B}}}
\put(2701,-1186){\makebox(0,0)[lb]{\smash{\SetFigFont{12}{14.4}{rm}$\alpha$}}}
\put(2521,434){\makebox(0,0)[lb]{\smash{\SetFigFont{12}{14.4}{rm}f}}}
\put(6031,-1051){\makebox(0,0)[lb]{\smash{\SetFigFont{12}{14.4}{rm}C}}}
\put(2476,-2716){\makebox(0,0)[lb]{\smash{\SetFigFont{12}{14.4}{rm}g}}}
\put(4771,-736){\makebox(0,0)[lb]{\smash{\SetFigFont{12}{14.4}{rm}h}}}
\end{picture}

gives a verical composite $ h \circ \alpha $. This is same as\\

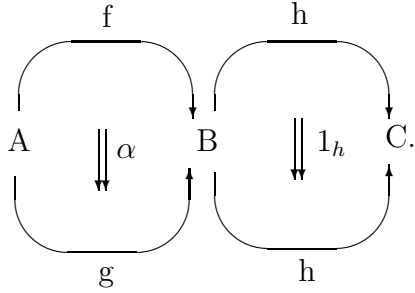
\setlength{\unitlength}{2000sp}%
\begingroup\makeatletter\ifx\SetFigFont\undefined
\def\x#1#2#3#4#5#6#7\relax{\def\x{#1#2#3#4#5#6}}%
\expandafter\x\fmtname xxxxxx\relax \def\y{splain}%
\ifx\x\y   
\gdef\SetFigFont#1#2#3{%
  \ifnum #1<17\tiny\else \ifnum #1<20\small\else
  \ifnum #1<24\normalsize\else \ifnum #1<29\large\else
  \ifnum #1<34\Large\else \ifnum #1<41\LARGE\else
     \huge\fi\fi\fi\fi\fi\fi
  \csname #3\endcsname}%
\else
\gdef\SetFigFont#1#2#3{\begingroup
  \count@#1\relax \ifnum 25<\count@\count@25\fi
  \def\x{\endgroup\@setsize\SetFigFont{#2pt}}%
  \expandafter\x
    \csname \romannumeral\the\count@ pt\expandafter\endcsname
    \csname @\romannumeral\the\count@ pt\endcsname
  \csname #3\endcsname}%
\fi
\fi\endgroup
\begin{picture}(4747,3330)(1351,-2716)
\thinlines
\put(3646,-691){\vector( 0,-1){0}}
\put(2558,-691){\oval(2176,1896)[tr]}
\put(2558,-601){\oval(2144,1716)[tl]}
\put(2517,-1411){\oval(2152,1890)[bl]}
\put(2517,-1321){\oval(2168,2070)[br]}
\put(3601,-1321){\vector( 0, 1){0}}
\put(4992,-1366){\oval(2152,1890)[bl]}
\put(4992,-1276){\oval(2168,2070)[br]}
\put(6076,-1276){\vector( 0, 1){0}}
\put(6076,-691){\vector( 0,-1){0}}
\put(4988,-691){\oval(2176,1896)[tr]}
\put(4988,-601){\oval(2144,1716)[tl]}
\put(2476,-826){\vector( 0,-1){765}}
\put(2566,-826){\vector( 0,-1){765}}
\put(4996,-691){\vector( 0,-1){765}}
\put(4906,-691){\vector( 0,-1){765}}
\put(1351,-1096){\makebox(0,0)[lb]{\smash{\SetFigFont{12}{14.4}{rm}A }}}
\put(3691,-1096){\makebox(0,0)[lb]{\smash{\SetFigFont{12}{14.4}{rm}B}}}
\put(2701,-1186){\makebox(0,0)[lb]{\smash{\SetFigFont{12}{14.4}{rm}$\alpha$}}}
\put(2521,434){\makebox(0,0)[lb]{\smash{\SetFigFont{12}{14.4}{rm}f}}}
\put(4861,479){\makebox(0,0)[lb]{\smash{\SetFigFont{12}{14.4}{rm}h}}}
\put(5176,-1096){\makebox(0,0)[lb]{\smash{\SetFigFont{12}{14.4}{rm}$1_{h}$}}}
\put(6031,-1051){\makebox(0,0)[lb]{\smash{\SetFigFont{12}{14.4}{rm}C.}}}
\put(2476,-2716){\makebox(0,0)[lb]{\smash{\SetFigFont{12}{14.4}{rm}g}}}
\put(4951,-2716){\makebox(0,0)[lb]{\smash{\SetFigFont{12}{14.4}{rm}h}}}
\end{picture}
\\

A {\em 2-functor} $F: \cal C \longrightarrow \cal D $ between
2-categories $\cal C$ and $\cal D$ is a triple of functions sending
0-cells, 1-cells, and 2-cells of $\cal C$ to items of the same types
in $\cal D$ so as to preserve all the categorical structures (domain,
codomain, identities, and composites).

A natural example of 2-categories is the category of {\bf 2-Top} whose
objects or 0-cells are topological spaces, 1-cells are continuous maps
between spaces, 2-cells are homotopy classes of homotopies between
continuous maps.

Another example is the category {\bf Grp} whose objects are groups,
1-cells are homomorphism between two groups, and 2-cells $\alpha : f
\Rightarrow g$ are the automorphims  of the codomain of $g$.

\section{semistrict monoidal 2-category}

(cf. [6], [7], [3])  A {\em semistrict monoidal 2-category} category\\ 

$\cal C$ = ( ${\cal C}$ , ${\otimes}$ , $I$, ${1 \otimes 0}$, ${0 \otimes
1}$,  ${2 \otimes 0}$, $ {0 \otimes 2}$, ${{\otimes}_{1,1'}} : {0 \otimes 1'}
\cong {0 \otimes 1'} $ whenever $(0 \otimes 1') (1 \otimes 0) = (1
\otimes 0) (0 \otimes 1')$ )\\
 consists of:\\

\begin{description}
\item[1.] A 2-category $\cal C$
\item[2.] For any two objects or 0-cells $A$ and $B$ in $\cal C$, an object
$A \otimes B$ in $\cal C$. 
\item[3.] The unit object $I \in \cal C$.
\item[4.] 1-cell composite functions
${1 \otimes 0}, {0 \otimes
1}$ and 2-cell composite functions ${2 \otimes 0}$,  ${0 \otimes 2}$
and $
{{\otimes}_{1,1'}} : {0 \otimes 1'}
\cong {0 \otimes 1'} $ whenever $(0 \otimes 1') (1 \otimes 0) = (1 \otimes 0) (0 \otimes 1')$, where 0,1 and 2 stand for the 0-cells,
1-cells and 2-cells respectively.
\end{description}
these composite functions act as follows :
\begin{description}
\item[4a.] For any 1-cell $f: A \longrightarrow B$ and any 0-cell $C \in \cal
C$ a 1-cell $f \otimes C : A \otimes B  \longrightarrow A \otimes C
$.
\item[4b.] For any 1-cell $g: B \longrightarrow C$ and any 0-cell $A \in \cal
C$ a 1-cell $A \otimes g : A \otimes B  \longrightarrow A \otimes C
$.
\item[4c.] For any 2-cell $\alpha : f \Rightarrow f'$ and any 0-cell $B \in \cal
C$ a 2-cell $\alpha \otimes B : f \otimes B  \Rightarrow f' \otimes B
$.
\item[4d.] For any 2-cell $\beta : g \Rightarrow g'$ and any 0-cell $A \in \cal
C$ a 2-cell $A \otimes \beta : A \otimes g  \Rightarrow A \otimes g'
$. 
\item[4e.] For any two 1-cells $f: A \longrightarrow B$ and $g: C
\longrightarrow D$  a  2-isomorphism  
${{\otimes}_{f,g}} : {A \otimes g}
\cong {B \otimes g} $\\

\setlength{\unitlength}{3500sp}%
\begingroup\makeatletter\ifx\SetFigFont\undefined
\def\x#1#2#3#4#5#6#7\relax{\def\x{#1#2#3#4#5#6}}%
\expandafter\x\fmtname xxxxxx\relax \def\y{splain}%
\ifx\x\y   
\gdef\SetFigFont#1#2#3{%
  \ifnum #1<17\tiny\else \ifnum #1<20\small\else
  \ifnum #1<24\normalsize\else \ifnum #1<29\large\else
  \ifnum #1<34\Large\else \ifnum #1<41\LARGE\else
     \huge\fi\fi\fi\fi\fi\fi
  \csname #3\endcsname}%
\else
\gdef\SetFigFont#1#2#3{\begingroup
  \count@#1\relax \ifnum 25<\count@\count@25\fi
  \def\x{\endgroup\@setsize\SetFigFont{#2pt}}%
  \expandafter\x
    \csname \romannumeral\the\count@ pt\expandafter\endcsname
    \csname @\romannumeral\the\count@ pt\endcsname
  \csname #3\endcsname}%
\fi
\fi\endgroup
\begin{picture}(3375,2718)(631,-2770)
\thinlines
\put(1756,-2266){\vector( 1, 0){1125}}
\put(1756,-601){\vector( 1, 0){1125}}
\put(3376,-871){\vector( 0,-1){900}}
\put(1441,-871){\vector( 0,-1){900}}
\put(2071,-1096){\vector( 0,-1){450}}
\put(2161,-1096){\vector( 0,-1){450}}
\put(1981,-196){\makebox(0,0)[lb]{\smash{\SetFigFont{12}{14.4}{rm}$A \otimes g$}}}
\put(991,-646){\makebox(0,0)[lb]{\smash{\SetFigFont{12}{14.4}{rm}$A \otimes C$}}}
\put(1036,-2356){\makebox(0,0)[lb]{\smash{\SetFigFont{12}{14.4}{rm}$B \otimes C$}}}
\put(4006,-1321){\makebox(0,0)[lb]{\smash{\SetFigFont{12}{14.4}{rm}$f \otimes D$}}}
\put(3601,-556){\makebox(0,0)[lb]{\smash{\SetFigFont{12}{14.4}{rm}$A \otimes D$}}}
\put(3556,-2356){\makebox(0,0)[lb]{\smash{\SetFigFont{12}{14.4}{rm}$B \otimes D$}}}
\put(631,-1366){\makebox(0,0)[lb]{\smash{\SetFigFont{12}{14.4}{rm}$f \otimes C$}}}
\put(2116,-2716){\makebox(0,0)[lb]{\smash{\SetFigFont{12}{14.4}{rm}$B \otimes g$}}}
\put(2476,-1411){\makebox(0,0)[lb]{\smash{\SetFigFont{12}{14.4}{rm}${\otimes}_{f,g}$.}}}
\end{picture}

\end{description}

such that the following conditions are satisfied\\

\begin{description}
\item[(i)] For any object $A \in \cal C$ we have {\em 2-functors} $A \otimes
  {-} : \cal C \longrightarrow \cal C$ and ${-} \otimes
  A : \cal C \longrightarrow \cal C$.
\item[(ii)] $A \otimes I = A = I \otimes A$ for any object $A$,\\
      $f \otimes I = f = I \otimes f$ for any 1-cell $f$,\\
$\alpha \otimes I = \alpha = I \otimes \alpha$ for any 2-cell $\alpha$,
\item[(iii)] Let $X$ be any object, 1-cell or 2-cell in $\cal c$, for
  all $X$, $A$ and $B$ in $\cal C$ we have\\
$(A \otimes B) \otimes X = A \otimes (B \otimes X)$, $(A \otimes X)
  \otimes B = A \otimes (X \otimes B)$,
$(X \otimes A) \otimes B = X \otimes (A \otimes B)$.
\item[(iv)] For any 1-cells $f : A \longrightarrow A'$, $g : B
  \longrightarrow B'$ and $h : C \longrightarrow C'$ in $\cal C$ we
  have
${\otimes}_{{A \otimes g}, h} =  A \otimes {\otimes}_{g, h}$,
  ${\otimes}_{{f \otimes B}, h} =  {\otimes}_{f, {B \otimes h}}$ and 
 ${\otimes}_{f, {g \otimes C}} =  {\otimes}_{f,g}  \otimes C $.
\item[(v)] For any objects $A$ and $B$ in $\cal C$ we have ${1_{A}} \otimes
  B = A \otimes {1_{B}} = 1_{A \otimes B} $.
\item[(vi)] For any 1-cells $f : A \longrightarrow A'$, any 1-cells $g, g' : B
  \longrightarrow B'$ and any 2-cell $\beta : g \Rightarrow g'$ we
  have 1-cell identities
\[(A' \otimes g) (f \otimes B) = (f \otimes B') (A \otimes g) \]
\[(A' \otimes g') (f \otimes B) = (f \otimes B') (A \otimes g') \]
and a 2-cell identity
\[(A' \otimes \beta) \circ ({\otimes}_ {f,g}) = ({\otimes}_{f,g'}) \circ (A \otimes \beta). \]
\item[(vii)]  For any 1-cells $g : B \longrightarrow B'$, any 1-cells $f, f' : A
  \longrightarrow A'$ and any 2-cell $\alpha : f \Rightarrow f'$ we
  have 1-cell identities
\[(A' \otimes g) (f' \otimes B) = (f' \otimes B') (A \otimes g) \]
\[(A' \otimes g) (f \otimes B) = (f \otimes B') (A \otimes g) \]
and a 2-cell identity
\[(\alpha \otimes B) \circ ({\otimes}_ {f,g}) = ({\otimes}_{f',g}) \circ (\alpha \otimes B). \]
\item[(viii)] For any 1-cells $f : A \longrightarrow A'$, $g : B
  \longrightarrow B'$ and $g' : B' \longrightarrow B''$ the
  2-isomorphism ${\otimes}_{f, gg'}$ equals to the pasting of
  ${\otimes}_{f, g}$ and ${\otimes}_{f, g'}$ as in the following diagram.\\

\setlength{\unitlength}{3500sp}%
\begingroup\makeatletter\ifx\SetFigFont\undefined
\def\x#1#2#3#4#5#6#7\relax{\def\x{#1#2#3#4#5#6}}%
\expandafter\x\fmtname xxxxxx\relax \def\y{splain}%
\ifx\x\y   
\gdef\SetFigFont#1#2#3{%
  \ifnum #1<17\tiny\else \ifnum #1<20\small\else
  \ifnum #1<24\normalsize\else \ifnum #1<29\large\else
  \ifnum #1<34\Large\else \ifnum #1<41\LARGE\else
     \huge\fi\fi\fi\fi\fi\fi
  \csname #3\endcsname}%
\else
\gdef\SetFigFont#1#2#3{\begingroup
  \count@#1\relax \ifnum 25<\count@\count@25\fi
  \def\x{\endgroup\@setsize\SetFigFont{#2pt}}%
  \expandafter\x
    \csname \romannumeral\the\count@ pt\expandafter\endcsname
    \csname @\romannumeral\the\count@ pt\endcsname
  \csname #3\endcsname}%
\fi
\fi\endgroup
\begin{picture}(5127,1890)(991,-2356)
\thinlines
\put(1756,-2266){\vector( 1, 0){1125}}
\put(1756,-601){\vector( 1, 0){1125}}
\put(3376,-871){\vector( 0,-1){900}}
\put(2071,-1096){\vector( 0,-1){450}}
\put(2161,-1096){\vector( 0,-1){450}}
\put(4501,-601){\vector( 1, 0){1125}}
\put(6076,-916){\vector( 0,-1){900}}
\put(4231,-2266){\vector( 1, 0){1125}}
\put(4861,-1141){\vector( 0,-1){450}}
\put(4771,-1141){\vector( 0,-1){450}}
\put(1261,-871){\vector( 0,-1){900}}
\put(991,-646){\makebox(0,0)[lb]{\smash{\SetFigFont{12}{14.4}{rm}$A \otimes B$}}}
\put(1036,-2356){\makebox(0,0)[lb]{\smash{\SetFigFont{12}{14.4}{rm}$A'\otimes
      B$}}}
\put(2476,-1411){\makebox(0,0)[lb]{\smash{\SetFigFont{12}{14.4}{rm}${\otimes}_{f,
          g}$}}}
\put(6031,-646){\makebox(0,0)[lb]{\smash{\SetFigFont{12}{14.4}{rm}$A
      \otimes B''$}}}
\put(5896,-2311){\makebox(0,0)[lb]{\smash{\SetFigFont{12}{14.4}{rm}$A'
      \otimes B''$}}}
\put(4996,-1366){\makebox(0,0)[lb]{\smash{\SetFigFont{12}{14.4}{rm}${\otimes}_{f,
          g'}$}}}
\put(3421,-601){\makebox(0,0)[lb]{\smash{\SetFigFont{12}{14.4}{rm}$A
      \otimes B'$}}}
\put(3331,-2356){\makebox(0,0)[lb]{\smash{\SetFigFont{12}{14.4}{rm}$A'
      \otimes B'$}}}
\end{picture}

\end{description}

Similarly, the 2-isomorphism ${\otimes}_{ff', g}$ equals to the pasting
of ${\otimes}_{f, g}$ and ${\otimes}_{f', g}$.

\section{Internal Categories}

We can define monoid, group, graph, and other structures in a category
$\cal C$. We can also define a category within $\cal C$ -called a {\em
  category of objects} in $\cal C$ or an {\em internal category} in
$\cal C$. Such an internal category provides a generalize  version of
the category $\cal C$.\\ 
In what follows we assume category $\cal C$ is finite complete.\\

An {\em internal category} {\bf C} in  $\cal C$ consists of : 
\begin{itemize}
\item an object $C_{0} \in  obj( {\cal C})$, called {\em object of objects};
\item an object of morphisms $C_{1} \in obj( {\cal C} )$, called {\em object of
  arrows};\\

together with four maps in $\cal C$
\item {\em source or domain morphism} $s: C_{1} \longrightarrow C_{0}$
  and {\em target or codomain morphism } $t: C_{1} \longrightarrow
  C_{0}$;
\item an {\em identity arrow}  $i: C_{0} \longrightarrow C_{1}$;
\item a {\em composition morphism} $\odot : {C_{1}} {{\times}_{C_{0}}}
  C_{1} \longrightarrow C_{1}$,
here composition $\odot$ is defined on the following pullback ($
  {C_{1}} {{\times}_{C_{0}}} C_{1}$, $p$, $q$) of $s$ and $t$ :

\setlength{\unitlength}{3500sp}%
\begingroup\makeatletter\ifx\SetFigFont\undefined
\def\x#1#2#3#4#5#6#7\relax{\def\x{#1#2#3#4#5#6}}%
\expandafter\x\fmtname xxxxxx\relax \def\y{splain}%
\ifx\x\y   
\gdef\SetFigFont#1#2#3{%
  \ifnum #1<17\tiny\else \ifnum #1<20\small\else
  \ifnum #1<24\normalsize\else \ifnum #1<29\large\else
  \ifnum #1<34\Large\else \ifnum #1<41\LARGE\else
     \huge\fi\fi\fi\fi\fi\fi
  \csname #3\endcsname}%
\else
\gdef\SetFigFont#1#2#3{\begingroup
  \count@#1\relax \ifnum 25<\count@\count@25\fi
  \def\x{\endgroup\@setsize\SetFigFont{#2pt}}%
  \expandafter\x
    \csname \romannumeral\the\count@ pt\expandafter\endcsname
    \csname @\romannumeral\the\count@ pt\endcsname
  \csname #3\endcsname}%
\fi
\fi\endgroup
\begin{picture}(2925,2403)(856,-2626)
\thinlines
\put(1756,-2266){\vector( 1, 0){1125}}
\put(1756,-601){\vector( 1, 0){1125}}
\put(3376,-871){\vector( 0,-1){900}}
\put(1261,-871){\vector( 0,-1){900}}
\put(3421,-601){\makebox(0,0)[lb]{\smash{\SetFigFont{12}{14.4}{rm} $C_{1}$}}}
\put(3331,-2356){\makebox(0,0)[lb]{\smash{\SetFigFont{12}{14.4}{rm}$
      C_{0}$ ;}}}
\put(901,-1276){\makebox(0,0)[lb]{\smash{\SetFigFont{12}{14.4}{rm}$p$}}}
\put(3781,-1366){\makebox(0,0)[lb]{\smash{\SetFigFont{12}{14.4}{rm}$t$}}}
\put(2206,-331){\makebox(0,0)[lb]{\smash{\SetFigFont{12}{14.4}{rm}$q$}}}
\put(856,-646){\makebox(0,0)[lb]{\smash{\SetFigFont{12}{14.4}{rm}${C_{1}} {{\times}_{C_{0}}} C_{1}$}}}
\put(901,-2356){\makebox(0,0)[lb]{\smash{\SetFigFont{12}{14.4}{rm} $C_{1}$}}}
\put(2206,-2626){\makebox(0,0)[lb]{\smash{\SetFigFont{12}{14.4}{rm}$s$}}}
\end{picture}

This is equal to the following two conditions :

$s \cdot \odot = s \cdot p$, and $t \cdot \odot = s \cdot q$.\\

These data must satisfy the following commutative conditions, which
simply express the usual axiom for a category :

\item $s$ $i$ = $1_{C_{0}}$ = $t$ $i $\\ 
specifies domain and codomain of the identity arrows;
\item $s \cdot \odot = s \cdot p$, and  $t \cdot \odot = t \cdot q$\\
  assigns the domain and codomain of composite morphisms;
\item $\odot$ $\circ$ ($\odot {{\times}_{C_{0}}} 1$) $ = $ $\odot$
  $\circ$ ($1 {{\times}_{C_{0}}} \odot$) \\
expresses that associative law for composition in terms of triple pullback;
\item $p =$ $\odot$ $\circ$ (${i \times 1} $), and  $ \odot$ $\circ$ (
  $ {1 \times i}$) $ = q $\\
gives the left and right unit laws for composition of morphisms.
\end{itemize}

When $\cal C = {\bf Set}$( category of small sets and functions) pullback is the set of composable pairs $(g,
f)$ of arrows.\\

An internal category in {\bf Set} is just an ordinary small
category which is same as an object in {\bf Cat} (category of small
categories and functors). 

An internal category in {\bf Grp} (category of groups and
homomorphisms)  is a category in which both $C_{0}$ and $C_{1}$ are
groups, and all the maps $i$, $s$, $t$ and $\odot$ are homomorphisms
of groups. We observe that internal category in {\bf Grp} is same as a
group object in {\bf Cat}.\\

An {\em internal functor}  (or {\em functor in} $\cal C$) $F : {\bf C}
\longrightarrow {\bf D}$ between two internal categoris {\bf C} and
{\bf D} of $\cal C$ consists of :\\
\begin{itemize}
\item  morphisms $F_{0} : C_{0} \longrightarrow D_{0}$, and $F_{1} :
  C_{1} \longrightarrow D_{1}$ of $\cal C$;\\
such that following holds :
\item $F_{0} \cdot s = s' \cdot F_{1} $, and $F_{0} \cdot t = t' \cdot
  F_{1} $\\
preservation of domain and codomain;
\item $F_{1} \cdot i = i' \cdot F_{0} $\\
preservation of identity arrows;
\item $F_{1} \circ {{\odot}_{\bf C}}  = {{\odot}_{\bf D}} \circ ( F_{1}
  \times F_{1}) $\\
preservation of composite arrows.
\end{itemize}

Using a similar procedure, an {\em internal natural
  transformation} between two internal functors $F$ and $G$from {\bf
  C} to {\bf D} in $\cal C$, say $\alpha : F \Rightarrow g$, is a
  morphism $\alpha : C_{0} \longrightarrow C_{1}$ which satisfies the
  following conditions :
\begin{itemize}
\item $s \cdot \alpha = F_{0}$, and $t \cdot \alpha = G_{0}$
\item $\odot \circ {\bigtriangleup (\alpha {\cdot} t {\times} F_{1} )} =
  \odot \circ {\bigtriangleup (G_{1} \times \alpha {\cdot} s )}$\\
where $\bigtriangleup$ is a diagonal morphism.
\end{itemize}

\subsection{2-category of internal categories}
\begin{pro} 
Internal categories, internal functors, and internal natural
transformations in $\cal C$ form a strict 2-category $ {\bf 2}{\cal C}$.
\end{pro}
We denote by {\bf Vect}, the category of vector spaces over a field $k$ and
linear functions.
\begin{defi}
A {\bf 2-vector space} is an internal category in {\bf Vect}.
\end{defi}

When $\cal C = $ {\bf Grp} or {\bf Vect}, we can formulate the
definition of an internal categories in $\cal C$ without the use of the
categorical composition $\odot$  in $\cal C$. (cf. [1]).

\begin{pro}
2-vector spaces, linear functors and linear natural transformations in
{\bf Vect} form a 2-category {\bf 2Vect} the 2-category of 2-vector
spaces .
\end{pro}

 Following the Gray monoid construction of Day and Street [5] we  prove 
\begin{pro}
{\bf 2Vect} forms a semistrict monoidal 2-category.
\end{pro}

Note that a discrete
2-vector space category is simply a vector space i.e.  an object of {\bf Vect}.

\section{Extended 2-dimensional TQFTs}

An $n$-dimensional TQFT is a symmetric monoidal functor $F : {\bf Cob}_{n+1}
\longrightarrow {\bf Vect}$. Here the category  ${\bf Cob}_{n+1}$ has
compact oriented $n$-dimensional manifolds as objects and compact
oriented cobordisms, which are equivalence classes of
$(n+1)$-manifolds with boundary, between them as
morphisms, and it has a monoidal structure (tensor product) given by 
disjoint union. An  example of physicists interest is the 2-category of
{\em relative cobordism}, ${{\bf Cod}^{rel}}_{n+1}$ (cf. [8]) has
the underlying category ${\bf Cob}_{n+1}$ as objects and 1-cells. The
2-cells between two cobordisms are given by $(n+2)$-dimensional
manifolds with boundary satisfying certain conditions. This category
can be formulated in a semistrict version of monoidal 2-category.
 
In $n$-categorical set up, other examples of monoidal 2-categories
are

1. category {\bf Chcomp} has chain complexes as 0-cells, 
chain maps as 1-cells, and chain homotopies as 2-cells.

2. category {\bf n-Cob} has 0-manifolds as 0-cells ( we assume all
manifolds are `compact, smooth, oriented manifolds'), 1-manifolds with
corners. i.e. cobordism between 0-manifolds as 1-cells, and 2-manifolds with
corners as 2-cells. (cf [9]).

Instead of taking 0-cells as 0-manifolds, one can also start with
objects as 1-manifolds with or without corners to get
Atiya-Segal-style TQFT.

A 2-dimensional TQFT is a particular case of the above
construction. Here the category ${\bf Cob}_{1+1}$ or ${\bf Cob}_{2}$
has compact oriented 1-manifolds as objects and compact oriented
cobordism between them as morphisms.\\

{\em Extended TQFTs} constructed by  Kerler and Lyubashenko in
[8]  involves higher category theory, namely double
categories and double functors. Their construction of extended version
of TQFTs is quite different from the $n$-categorical version of
extended TQFTs purposed by the Baez and  Dolan in [2].

Reshetikhin and Turaev in 1990
constructed a ribbon invariants defined via quantum groups and later in
1991 they succeeded for the first time to construct 3-manifold
invariants in a rigorous and mathematically consistent way. In their
paper they obtained a projective TQFT in the sense of Atiyah. They
use a {\em semisimple modular category}  as input data.
Turaev's  monograph [10] fully describe the Reshetikhin and Turaev
construction of a TQFT, methods used in [10] show that the construction
of a TQFT functor is
very complicated and a difficult procedure.\\

Extended TQFT functor
constructed in [8] (which is more complicated then  the ordinary TQFT
functor of [T]) cannot be considered as a generalized version of Turaev
construction of  TQFT functor, actually both construction are
different because of the different base categories. 

Baez and Dolan's [2]
hypothesis for  extended TFQTs shows that one needs a different
construction of  TQFT functors at a  different dimensional level,
e.g. the TQFT functor which produce 2-dimensional extended  TQFTs
cannot be (easily) generalized to a 3-dimensional extended
TQFT functor. 
This is because of the base categories (especially
the codomain categories). Either they do not have a nice structure in higher
dimension or their structure is very complicated, e.g. the enriched
$n$-categorical version of {\bf Vect} is not very clear in 
dimension $n \geq 2$. 

For $n = 2$ , one can think  {\em 2-vector spaces} ( in the sense of
$n$-categories)  as
a linear space over the category {\bf Vect} of vector spaces. This means
that, it is a monoidal category $\cal V$ with an external tensor product
$\oplus$ and a functor $\otimes : {\bf Vect}  \times \cal V
\longrightarrow \cal V$ satisfying various conditions. Involvement of
external tensor product over categories suggests that one needs to construct different
TQFT functors at different dimensional level. This also suggests that in
most of the higher dimensional cases these TQFTs functors will be independent from each other.

In a general situation one can ask ``is there any other way to construct an
$n$-dimensional extended TQFTs such   that 
\begin{itemize}
\item an $(n-1)$-dimension extended TQFT can be obtained by putting
  certain constraints on the $n$-dimensional extended TQFT functor?
\item a generalized version of $(n-1)$-dimensional
  extended TQFT can be realized.
\end{itemize}
As an attempt to answer this question, we define 2-dimensional
extended TQFTs in the following hypothesis:
\begin{defi}
A two dimensional extended TQFT is a functor from semistrict monoidal
2-category of {\bf 2-cob} to semistrict monoidal 2-category {\bf 2Vect}
of 2-vector spaces. 
\end{defi}
Here, internal categories {\em 2-Vector spaces} are same as  the
2-Vector spaces defined  by Baez and Crans in [1].\\

For the higher dimensional extended TQFTs, one needs 
to generalize internal categories structure for higher
dimensions in such a way that existing base category structures remain
preserved, e.g.
as in the case of {\bf 2Vect}, which contains ordinary vector spaces
as objects. If we consider {\bf 3Vect} to be the category having
objects as internal categories of {\bf 2Vect} and arrows are  internal
functors, then under suitable conditions {\bf 3Vect} can gives a higher
category version of {\bf 2Vect} which also contain {\em 2-vector
  spaces}  as objects.
 
Presence of an ordinary vector space in {\bf 2Vect} is vital for
following two reasons 

1. to get a Turaev's  type of  TQFT from
a 2-dimensional extended TQFT.\\

2. to generalize Turaev's type TQFT to 2-dimensional extended
TQFT.\\ 
These  two cases can be proved by making  a
suitable modification or restriction in TQFT functors. 

One can also think an {\em n-categories} version of this result in
terms of enriched categories by using the Baez's construction of
$n$-categories.

Details on the various  aspects of  2-dimensional extended TQFT which
are discussed here is a subject matter of the next  version of this article.


\begin{thebibliography}{1}
\bibitem{bc}J. Baez and A.S. Crans,
Higher-dimensional algebra VI: Lie 2-algebra, math:QA/0307263,Jul
2003.
\bibitem{bd}J. Baez and J. Dolan, Higher-dimensional algebra and topological quantum field theory.
J. Math. Phys.  36  (1995),  no. 11, 6073--6105.
\bibitem{bn} J.  Baez and M.  Neuchl, Higher-dimensional
  algebra. I. Braided monoidal $2$-categories.  Adv. Math.  121
  (1996),  no. 2, 196--244.
\bibitem{C} L. Crane, Clock and category: is quantum
  gravity algebraic?  J. Math. Phys.  36  (1995),  no. 11,
  6180--6193. 
\bibitem{ds} B. Day and  R. Street, Monoidal bicategories and Hopf algebroids.  Adv. Math.  129  (1997),  no. 1, 99--157. 
\bibitem{gs} R. Gordon, A. J. Power and R.  Street, Coherence for tricategories.  Mem. Amer. Math. Soc.  117  (1995),  no. 558, vi+81 pp.
\bibitem{kv} M. M. Kapranov and V. A.  Voevodsky, $2$-categories and Zamolodchikov tetrahedra equations.  Algebraic groups and their generalizations: quantum and infinite-dimensional methods (University Park, PA, 1991),  177--259, Proc. Sympos. Pure Math., 56, Part 2, Amer. Math. Soc., Providence, RI, 1994. 
\bibitem{tl}T. Kerler and
V.V. Lyubashenko, Non-semisimple topological quantum field theories
for 3-manifolds with corners. Lecture Notes in Mathematics,
1765. Springer-Verlag, Berlin, 2001.
\bibitem{Li} T. Leinster, Topology and Higher-dimensional category theory:
  the rough idea, math:CT/0106240, Jun 2001.
\bibitem{T} V. G. Turaev, Quantum invariants of knots and 3-manifolds. de Gruyter Studies in Mathematics, 18. Walter de Gruyter  Co., Berlin, 1994. x+588 pp. 

\end{thebibliography}
\end{document}